\documentclass[11pt,twoside]{article}
\usepackage{mathrsfs}
\usepackage{amsmath}
\usepackage{amssymb}
\usepackage{graphicx,color}
\usepackage{subfigure}
\usepackage{appendix,url}
\newtheorem{theorem}{Theorem}

\newtheorem{proposition}{Proposition}
\newtheorem{claim}{Claim}
\newtheorem{remark}{Remark}

\numberwithin{equation}{section}

\title{Evolution equations on time-dependent intervals}
\author{A.S. Fokas$^{1}$, B. Pelloni$^{2}$ and Baoqiang Xia$^{3}$
\\
$^1$Department of Applied Mathematics and Theoretical Physics,
\\University of Cambridge, Cambridge CB3 0WA, United Kingdom,
\\
 E-mail address:
t.fokas@damtp.cam.ac.uk
\\
$^2$Department of Mathematics, Heriot-Watt University,
\\Edinburgh EH14 4AS, United Kingdom,
\\
 E-mail address:
 b.pelloni@hw.ac.uk
 \\
 $^{3}$School of Mathematics and Statistics, Jiangsu Normal
University,\\
 Xuzhou, Jiangsu 221116, P. R. China,
 \\
 E-mail address:
xiabaoqiang@126.com}

\date{\today}
\begin{document}
\maketitle

\begin{abstract}
{We study initial boundary value problems for linear evolution partial differential equations (PDEs)
posed on a time-dependent interval $l_1(t)<x<l_2(t)$, $0<t<T$,
where $l_1(t)$ and $l_2(t)$ are given, real, differentiable functions, and $T$ is an arbitrary constant.
For such problems, we show how to characterise the unknown boundary values in terms of the given initial and boundary conditions.
As illustrative examples we consider the heat equation and the linear Schr\"{o}dinger equation. In the first case, the unknown Neumann boundary values are expressed in terms of the Dirichlet boundary values and of the initial value through the unique solution of a system of two linear  integral equations with explicit kernels. In the second case, a similar result can  be proved  but only for a more restrictive class of boundary curves.} 
{Linear evolution PDE, Initial-boundary value problem, Riemann-Hilbert problem.}

%\noindent{\bf PACS:}\quad 02.30.Ik, 04.20.Jb.
\end{abstract}

\section{ Introduction}

We study linear evolution PDEs posed on a time-dependent interval. Specifically, we consider linear evolution PDEs on a $t$-dependent domain $\Omega(t)$ of the form
\begin{eqnarray}
\Omega(t)=\left\{(x,s):~l_1(s)<x<l_2(s),\;0<s<t\right\}\subset\mathbb{R}^2 \label{MI}
\end{eqnarray}
where $\left\{l_j(t)\right\}_1^2$ are given, real, continuously differentiable functions, such that  $l_1(s)<l_2(s)$ for all $s>0$ and $l_1(0)=0$, $l_2(0)=L\geq 0$. 

We first present some results for a general linear evolution PDE, and then concentrate on two illustrative examples, namely
the heat equation and the  linear Schr\"{o}dinger (LS) equation:
\begin{eqnarray}
(heat)~~q_t-q_{xx}=0,
\label{heat}
\\
(LS)~~iq_t+q_{xx}=0.
\label{ls}
\end{eqnarray}

For a given constant  $T>0$, we consider the above PDEs in the domain $\Omega(T)$ and assume that the following initial and Dirichlet boundary conditions are prescribed:
\begin{eqnarray}
&&q(x,0)=q_0(x)\in{\bf C}^1[0,L], ~~0<x<L,\label{iv}
\\
&&q(l_1(t),t)=f_0(t)\in {\bf C}^1[0,T],~~q(l_2(t),t)=g_0(t)\in {\bf C}^1[0,T], ~~0<t<T.
\label{bvD}
\end{eqnarray}
%%To ensure the unique existence of a solution, we assume the following regularity for the given data:
%\begin{itemize}
%\item[(a)]
%$q_0(x)\in{\bf C}^1[0,L]$;
%\item[(b)]
%$f_0(t), \,g_0(t)\in {\bf C}^1[0,T]$.
%\end{itemize}
We will develop our analysis in this general case, but for the particular example of the LS equation, we will need to restrict the class of boundary functions $l_1(t)$, $l_2(t)$. 

%The integral representation of $q(x,t)$, in addition to involving the above known functions,
%it also involves the unknown Neumann boundary values $q_x(l_1(t),t)$ and $q_x(l_2(t),t)$.
To obtain an effective representation of the solution  $q(x,t)$,  for $t<T$, one needs to 
determine the unknown Neumann boundary values $\left\{q_x(l_j(t),t)\right\}_1^2$, in terms of the given initial and Dirichlet boundary  data $\{q_0(x),f_0(t),g_0(t)\}$, i.e. to characterise the so-called Dirichlet to Neumann map.  The discussion and results in \cite{Xia} imply  that once this characterisation is achieved,  the solution obtained is in fact the {\em unique} solution of the boundary value problem.

The main result of this paper is the characterisation of the Dirichlet to Neumann map
%It is shown in sections 3 and 4 respectively, that
 for equations (\ref{heat}) and (\ref{ls})
% . Indeed we show, in sections 3 and 4 respectively, that $\left\{q_x(l_j(t),t)\right\}_1^2$
%can be obtained as 
through the unique solution of a system of two linear Volterra  integral equations.

The results presented here provide a generalization of the results of \cite{AP3} (see also \cite{DF}, \cite{FP} and \cite{P}), where
the analogous problem formulated on $l(t)<x<\infty$ was analysed. There are two important differences between that case and the results presented here. Firstly, we must characterise two, rather than one, unknown boundary values, and hence even in the case of second order linear equations the solution is given in terms of the solution of a {\em system} of Volterra linear integral equations.   Secondly,  the kernel of the  integral equations that characterise the unknown Neumann boundary value may be strongly singular. Hence  to obtain a rigorous existence result this kernel must  be regularised. 
%In addition to the analysis of equations (\ref{ls}) and (\ref{heat}), \cite{AP3} also contains the analysis of the linear KdV equation
%\begin{eqnarray}
%q_t+q_{xxx}=0.
%\label{kdv}
%\end{eqnarray}
%It was shown in \cite{AP3} that the analysis of (\ref{kdv}) is conceptually different with that the analysis of (\ref{ls}) and (\ref{heat}):
%in the latter case the relevant linear integral equations can be obtained directly by employing the simple integral representation  for $q(x,t)$
%obtained via the Fourier transform. On the other hand, in the case of (\ref{kdv}) one has to exploit the invariant properties of the so-called global relation.
%This relation plays a crucial role in the general methodology introduced in \cite{F1} and \cite{F2000} for solving boundary value problems known as the
%\emph{unified transform} or the \emph{Fokas method}.
%
%For equation (\ref{kdv}) formulated in $\Omega(T)$, the analysis of the associated global relation yield a system of four linear singular integral
%equations for the unknown functions  $\left\{q_x(l_j(t),t), q_{xx}(l_j(t),t)\right\}_1^2$, in terms of $\{q_0(x),f_0(t),g_0(t)\}$;
%these equations will be presented elsewhere.

The main results of this paper for the two equations (\ref{heat}) and (\ref{ls})  are Theorem \ref{main} and Theorem \ref{main2}, which use the following formal representation results:
\begin{proposition}[Heat equation]\label{theoremheat}
Let $q(x,t)$ be the solution of the heat equation (\ref{heat}) satisfying the initial
and Dirichlet boundary conditions (\ref{iv}) and (\ref{bvD}).
Denote by $f_1(t)$ and $g_1(t)$  the unknown Neumann boundary values $q_x(x,t)$ evaluated at $x=l_1(t)$ and at $x=l_2(t)$:
\begin{eqnarray}
f_1(t)=q_x(l_1(t),t), ~~g_1(t)=q_x(l_2(t),t),~~0<t<T.
\label{bvN}
\end{eqnarray}
The functions $f_1(t)$ and $g_1(t)$ can be expressed in
terms of the given initial and boundary data as the solution of the following system of coupled linear  integral equations:
\begin{subequations}
\begin{align}
&\pi f_1(t)=N_1(t)+\int_{0}^{t}K_{11}(t,s)f_1(s)ds-\int_0^tK_{12}(t,s)g_1(s)ds,~~0<t<T,
\label{qxrepresentationl1Fa}
\\
&\pi g_1(t)=N_2(t)-\int_{0}^{t}K_{22}(t,s)g_1(s)ds+\int_0^tK_{21}(t,s)f_1(s)ds,~~0<t<T,
\label{qxrepresentationl1Fb}
\end{align}
\label{qxrepresentationl1F}
\end{subequations}
%%K_{21}=-K_{12}
where the known functions $N_j(t)$, $j=1,2$, are given by
\begin{eqnarray}
N_j(t)=\sqrt{\pi}\Bigg\{\frac{1}{\sqrt{t}}\int_{0}^{L}e^{-\frac{(x-l_j(t))^2}{4t}}q'_0(x)dx
-\int_{0}^{t}\Big[\frac{e^{-\frac{(l_j(t)-l_1(s))^2}{4(t-s)}}}{\sqrt{t-s}}f_0'(s)
-\frac{e^{-\frac{(l_j(t)-l_2(s))^2}{4(t-s)}}}{\sqrt{t-s}}g_0'(s)\Big]ds\Bigg\}, 
\label{N1heat}
\end{eqnarray}
and the kernels $K_{jm}(t,s)$, $j,m=1,2$, are given by
\begin{equation}
K_{jm}(t,s)=\frac{\sqrt{\pi}}{2}\frac{l_j(t)-l_m(s)}{t-s}
\frac{e^{-\frac{(l_j(t)-l_m(s))^2}{4(t-s)}}}{\sqrt{t-s}},
~~0<s<t<T, ~~j=1,2.
\label{Kheat}
\end{equation}
%while the kernels $K_{12}(t,s,\varepsilon)$, $K_{21}(t,s,\varepsilon)$  are given by
%\begin{eqnarray}
%\begin{split}
%&K_{12}(t,s,\varepsilon)=\frac{\sqrt{\pi}}{2}\frac{l_1(t)-l_2(s)}{t-s+i\varepsilon}
%\frac{e^{-\frac{(l_1(t)-l_2(s))^2}{4(t-s+i\varepsilon)}}}{\sqrt{t-s+i\varepsilon}},
%~~\varepsilon>0,~~0<s<t<T,
%\end{split}
%\label{K1heat}
%\\
%\begin{split}
%&K_{21}(t,s,\varepsilon)=\frac{\sqrt{\pi}}{2}\frac{l_2(t)-l_1(s)}{t-s+i\varepsilon}
%\frac{e^{-\frac{(l_2(t)-l_1(s))^2}{4(t-s+i\varepsilon)}}}{\sqrt{t-s+i\varepsilon}},
%~~\varepsilon>0,~~0<s<t<T,.
%\end{split}
%\label{K2heat}
%\end{eqnarray}
%The kernels $K_{jm}$ are either regular (for $j\neq m$) or weakly singular (for $j=m$). Hence the system of Volterra integral equations (\ref{qxrepresentationl1F}) admits a unique solution $(f_1(t), g_1(t))$. 
\end{proposition}

\begin{proposition}[Linear Schr\"odinger equation]\label{theoremLS}

Let $q(x,t)$ be the solution of the linear Schr\"{o}dinger equation (\ref{ls}) satisfying the initial and Dirichlet boundary conditions (\ref{iv}) and (\ref{bvD}).
Let $f_1(t)$ and $g_1(t)$ denote the unknown Neumann boundary values as given in expression (\ref{bvN}).
The unknown boundary values $f_1(t)$ and $g_1(t)$ can be expressed in
terms of the given initial and boundary data as the solution of the following system of coupled linear  integral equations:
\begin{subequations}
\begin{align}
&\pi f_1(t)=N_1(t)+\int_{0}^{t}K_{11}(t,s)f_1(s)ds-\lim_{\varepsilon\to 0}\int_0^tK_{12}(t,s,\varepsilon)g_1(s)ds,~~0<t<T,
\label{qxrepresentationls1Fa}
\\
&\pi g_1(t)=N_2(t)+\lim_{\varepsilon\to 0}\int_{0}^{t}K_{21}(t,s,\varepsilon)f_1(s)ds-\int_0^tK_{22}(t,s)g_1(s)ds,~~0<t<T,
\label{qxrepresentationls1Fb}
\end{align}
\label{qxrepresentationls1F}
\end{subequations}
where the functions $N_j(t)$, $j=1,2$, are given by 
\begin{eqnarray}
N_j(t)=\frac{(1-i)\sqrt{2\pi}}{2}\Bigg\{\frac{1}{\sqrt{t}}\int_{0}^{L}e^{\frac{i(x-l_j(t))^2}{4t}}q'_0(x)dx
-\int_{0}^{t}\Big[\frac{e^{\frac{i(l_j(t)-l_1(s))^2}{4(t-s)}}}{\sqrt{t-s}}f_0'(s)
-\frac{e^{\frac{i(l_j(t)-l_2(s))^2}{4(t-s)}}}{\sqrt{t-s}}g_0'(s)\Big]ds\Bigg\}.
\nonumber \\
\label{N1ls}
\end{eqnarray}
%\end{claim} 
The singular  integral kernels $K_{jj}(t,s)$ are given by 
\begin{equation}
K_{jj}(t,s)=\frac{(1-i)\sqrt{2\pi}}{4}\frac{l_j(t)-l_j(s)}{t-s}
\frac{e^{\frac{i(l_j(t)-l_j(s))^2}{4(t-s)}}}{\sqrt{t-s}},
~~\varepsilon>0,~~0<s<t<T,
\label{K1ls}
\end{equation}
and the integral kernels $K_{jm}(t,s,\varepsilon)$, $j,m=1,2$, $j\neq m$ are given by
\begin{equation}
K_{jm}(t,s,\varepsilon)=\frac{(1-i)\sqrt{2\pi}}{4}\frac{l_j(t)-l_m(s)}{t-s-i\varepsilon}
\frac{e^{\frac{i(l_j(t)-l_m(s))^2}{4(t-s-i\varepsilon)}}}{\sqrt{t-s-i\varepsilon}},
~~\varepsilon>0,~~0<s<t<T.
\label{K2ls}
\end{equation} 
%The kernels $K_{jm}$ are either regular (for $j\neq m$) or weakly singular (for $j=m$). Hence the system of Volterra integral equations (\ref{qxrepresentationl1F}) admits a unique solution $(f_1(t), g_1(t))$. 
\end{proposition}
  
The representations above are formal, and it is not immediately clear that they actually define any function, let alone the boundary values of the given problem. In particular, in Proposition \ref{theoremLS}, we need to guarantee that the limits as $\varepsilon\to 0$ exist, and that the Volterra integral equations are well posed and admit a unique solution. 

Our main result for the heat equation is obtained by invoking classical theory for Volterra integral equations. 
\begin{theorem}\label{main}
For given functions $q_0(x)$, $f_0(t)$ and $g_0(t)$ as in (\ref{iv})-(\ref{bvD}),  the system of Volterra integral equations (\ref{qxrepresentationl1F})   admits a unique solution $(f_1(t), g_1(t))\in\,{\bf C}^1[0,T)\times {\bf C}^1[0,T)$. 
\end{theorem}

For the LS equation, we need to prove that the limits as $\varepsilon\to 0$ appearing (\ref{qxrepresentationls1F}) yield a regular kernel, and only then it will be possible to invoke classical theory.

\begin{theorem}\label{main2}
Assume  the the boundary functions $l_1(t)$, $l_2(t)$ are twice differentiable in $[0,T]$, and satisfy
\begin{equation}
%l_1(t)=\alpha t, \quad l_2(t)=\beta t+L,\qquad L>0,\;\;0<2\alpha<\beta.
%\label{lines}
l_1'(t)<0,\;\;l_1''(t)\geq 0;\qquad l_2'(t)>0,\;\;l_2''(t)\leq 0.
\label{bcurves}
\end{equation}
For given functions $q_0(x)$, $f_0(t)$ and $g_0(t)$ as in (\ref{iv})-(\ref{bvD}),  the limit system (\ref{qxrepresentationls1F}) is a system of Volterra integral equations (\ref{qxrepresentationls1F})  which admits a unique solution $(f_1(t), g_1(t))\in\,{\bf C}^1[0,T)\times {\bf C}^1[0,T)$. 
\end{theorem}

The paper is organized as follows.
In section 2, we derive formally a representation of $q(x,t)$ for a general evolution PDE formulated in $\Omega(T)$, as well as the associated global relation.
In section 3, we  derive the formal representations of Propositions  \ref{theoremheat} and \ref{theoremLS}. Finally, in section 4 we prove Theorem \ref{main} and Theorem \ref{main2}. In the case of the linear Schr\"odinger equation, we also remark on the important special case of linear boundaries. 

%The present paper is supplemented by paper \cite{Xia} where it is shown that if there exist functions $\{q_0(x),f_0(t),g_0(t),f_1(t),g_1(t)\}$ which satisfy
%the global relation, then $q(x,t)$ exists and furthermore
%\begin{eqnarray}
%\begin{split}
%&q(x,0)=q_0(x),~~q(l_1(t),t)=f_0(t),~~q(l_2(t),t)=g_0(t),
%\\
%&\partial_x q(l_1(t),t)=f_1(t),~~\partial_xq(l_2(t),t)=g_1(t).
%\end{split}
%\label{ibvls}
%\end{eqnarray}

\section{A formal integral representation for a general evolution PDE}
We consider the general linear evolution PDE
\begin{eqnarray}
&\left(\partial_t+i\sum_{j=1}^{n}\alpha_j(-i\partial_x)^j\right)q(x,t)=0, ~~(x,t)\in \Omega(T),
\label{EQ}
\end{eqnarray}
where $\alpha_n\neq 0$ and all $\alpha_j$'s are constants, and the domain $\Omega(T)$ is described by (\ref{MI}).

Let
\begin{eqnarray}
&&\omega(\lambda)=\sum_{j=1}^{n}\alpha_j \lambda^j, %~~\alpha_j\in \mathbb{R},
\label{dr}
\\
&&Q(x,t,\lambda)=-\sum_{j=1}^{n}\alpha_j\left((-i\partial_x)^{j-1}+\lambda(-i\partial_x)^{j-2}+\cdots+\lambda^{j-1}\right)q(x,t).
\label{Q}
\end{eqnarray}
The PDE (\ref{EQ}) can be written in the following divergence form:
\begin{eqnarray}
\left(e^{-i\lambda x+i\omega(\lambda )t}q(x,t)\right)_t=\left(e^{-i\lambda x+i\omega(\lambda )t}Q(x,t,\lambda )\right)_x.
\label{div}
\end{eqnarray}
Using the two-dimensional Green's theorem in the domain $\Omega(t)$, we obtain
\begin{eqnarray}
\oint_{\partial \Omega(t)} \left(e^{-i\lambda x+i\omega(\lambda )s}q(x,s)dx+e^{-i\lambda x+i\omega(\lambda )s}Q(x,s,\lambda )ds\right)=0,~~0<t<T,
\label{GRt}
\end{eqnarray}
where $\partial \Omega(t)$ denotes the oriented boundary of the domain $\Omega(t)$, such that $\Omega(t)$ lies to the left-hand side of the increasing direction.
Equation (\ref{GRt}) yields the relation
\begin{eqnarray}
&\int_{0}^{L}e^{-i\lambda x}q(x,0)dx-e^{i\omega(\lambda )t}\int_{l_1(t)}^{l_2(t)}e^{-i\lambda x}q(x,t)dx
%\nonumber \\&
-\int_{0}^{t}e^{-i\lambda l_1(s)+i\omega(\lambda )s}\left(q(l_1(s),s)l'_1(s)+Q(l_1(s),s,\lambda )\right)ds
\nonumber \\
&+\int_{0}^{t}e^{-i\lambda l_2(s)+i\omega(\lambda )s}\left(q(l_2(s),s)l'_2(s)+Q(l_2(s),s,\lambda )\right)ds=0, ~~\lambda \in \mathbb{C},~~0<t<T.
\label{GRt1}
\end{eqnarray}

Let
\begin{subequations}
\begin{align}
&\hat{q}_0(\lambda )=\int_{0}^{L}e^{-i\lambda x}q(x,0)dx,
\label{hatq0}
\\
&\hat{q}(t,\lambda )=\int_{l_1(t)}^{l_2(t)}e^{-i\lambda x}q(x,t)dx,
\label{hatqt}
\\
&\hat{Q}_1(t,\lambda )=\int_{0}^{t}e^{-i\lambda l_1(s)+i\omega(\lambda )s}\left(q(l_1(s),s)l'_1(s)+Q(l_1(s),s,\lambda )\right)ds,
\label{hatQ1t}
\\
&\hat{Q}_2(t,\lambda )=\int_{0}^{t}e^{-i\lambda l_2(s)+i\omega(\lambda )s}\left(q(l_2(s),s)l'_2(s)+Q(l_2(s),s,\lambda )\right)ds.
\label{hatQ2t}
\end{align}
\end{subequations}
Equation (\ref{GRt1}) can be rewritten in the form of the following global relation:
\begin{equation}
\hat{q}(t,\lambda )=e^{-i\omega(\lambda )t}\hat{q}_0(\lambda )-e^{-i\omega(\lambda )t}\hat{Q}_1(t,\lambda )+e^{-i\omega(\lambda )t}\hat{Q}_2(t,\lambda ), ~~\lambda \in \mathbb{C},~~0<t<T.
\label{GRt2}
\end{equation}

Equation (\ref{GRt2}) can be viewed either as the formal representation of the
solution, or as the starting point for determining the unknown boundary values.
Indeed, the term $\hat{q}(t,\lambda )$ in (\ref{GRt2})
is the Fourier transform of $q(x,t)$ on the finite interval $l_1(t)<x<l_2(t)$.
Inverting this Fourier transform for $q(x, t)$, we obtain the following
formal representation of the solution:
\begin{equation}
q(x,t)=\frac{1}{2\pi }\int_{-\infty}^{\infty}e^{i\lambda  x-i\omega(\lambda )t}\Big[\hat{q}_0(\lambda )-\hat{Q}_1(t,\lambda )+\hat{Q}_2(t,\lambda )\Big]d\lambda , ~~(x,t)\in\Omega(T).
\label{DFSolution}
\end{equation}

Assuming that $q(x,t) = 0$ for $x < l_1(t)$ and for $x > l_2(t)$, equation (\ref{DFSolution}) is also formally valid at
$x = l_1(t)$ and at $x = l_2(t )$. Hence evaluating the inverse Fourier transform (\ref{DFSolution}) at these two points, we obtain
\begin{equation}
q(l_1(t),t)=\frac{1}{\pi }\int_{-\infty}^{\infty}e^{i\lambda  l_1(t)-i\omega(\lambda )t}\Big[\hat{q}_0(\lambda )-\hat{Q}_1(t,\lambda )+\hat{Q}_2(t,\lambda )\Big]d\lambda ,
\label{DFSolutiong0}
\end{equation}
and
\begin{equation}
q(l_2(t),t)=\frac{1}{\pi }\int_{-\infty}^{\infty}e^{i\lambda  l_2(t)-i\omega(\lambda )t}\Big[\hat{q}_0(\lambda )-\hat{Q}_1(t,\lambda )+\hat{Q}_2(t,\lambda )\Big]d\lambda .
\label{DFSolutionf0}
\end{equation}

We note that in the paper \cite{Xia}  it is shown that if there exist sufficiently regular functions $$\{q_0(x),f_0(t),g_0(t),f_1(t),g_1(t)\}$$ which satisfy
the global relation, then there exists a unique regular solution $q(x,t)$ of the PDE such that 
\begin{eqnarray}
&q(x,0)=q_0(x),~~q(l_1(t),t)=f_0(t),~~q(l_2(t),t)=g_0(t),
\nonumber \\
&\partial_x q(l_1(t),t)=f_1(t),~~\partial_xq(l_2(t),t)=g_1(t).
\label{ibvls}
\end{eqnarray}

Therefore if the Dirichlet to Neumann map is constructed starting from the  assumption that the global relation holds, its solution does indeed provide the unique solution of the boundary value problem. 

\section{The integral equations - formal derivation}
\subsection{The heat equation}

We consider the heat equation (\ref{heat}) formulated in the time-dependent domain (\ref{MI}),
with the given initial value (\ref{iv}) and Dirichlet boundary conditions (\ref{bvD}).
The functions $f_1(t)$ and $g_1(t)$, as in (\ref{bvN}),  denote the unknown Neumann boundary values at $x=l_1(t)$ and $x=l_2(t)$ respectively.

Our aim is to 
determine the unknown boundary values $f_1(t)$ and $g_1(t)$ in terms of the given functions $q_0(x)$, $f_0(t)$ and $g_0(t)$, hence to characterise the Dirichlet to Neumann map.
In order to determine this map, we {\em solve the global relation for the unknown boundary values}, $f_1(t)$ and $g_1(t)$.

In the case of heat equation, in the notation of the previous section we have
\begin{eqnarray}
\omega(\lambda )=-i\lambda ^2, ~~Q(x,t)=q_x(x,t)+i\lambda q(x,t).
\label{heatOQ}
\end{eqnarray}
Then the global relation (\ref{GRt2}) becomes
\begin{eqnarray}
&\hat{q}_0(\lambda )-e^{\lambda ^2t}\hat{q}(t,\lambda )-\int_{0}^{t}e^{-i\lambda l_1(s)+\lambda ^2s}\left[\Big(i\lambda+l'_1(s)\Big)q(l_1(s),s)+q_x(l_1(s),s)\right]ds \nonumber
\\&+\int_{0}^{t}e^{-i\lambda l_2(s)+\lambda ^2s}\left[\Big(i\lambda+l'_2(s) \Big)q(l_2(s),s)+q_x(l_2(s),s)\right]ds=0, ~~\lambda \in \mathbb{C},~~0<t<T,
\label{GRt2heat}
\end{eqnarray}
where $\hat{q}_0(\lambda )$ and $\hat{q}(t,\lambda )$ are defined by (\ref{hatq0}) and (\ref{hatqt}).

In order to derive a representation of $q_x (x, t)$, we follow \cite{AP3} and first multiply (\ref{GRt2heat})
by $i\lambda $. For the terms involving the unknown boundary functions, integration by parts yields the following identity, with $j=1,2$:
%\begin{equation}
\begin{eqnarray}
%\begin{split}
%&
\int_{0}^{t}e^{-i\lambda l_j(s)+\lambda ^2s}\Big(-\lambda^2+i\lambda l'_j(s)\Big)q(l_j(s),s)ds
%\\
&=&e^{-i\lambda l_j(0)}q(l_j(0),0)-e^{\lambda^2t-i\lambda l_j(t)}q(l_j(t),t)
\nonumber \\
&&+\int_{0}^{t}e^{\lambda ^2s-i\lambda l_j(s)}\frac {dq}{ds}(l_j(s),s)ds.
%\end{split}
\label{dheat1}
\end{eqnarray}
%and
%\begin{eqnarray}
%\begin{split}
%&\int_{0}^{t}e^{-i\lambda l_2(s)+\lambda ^2s}\Big(-\lambda^2+i\lambda l'_2(s)\Big)q(l_2(s),s)ds
%\\
%=&q(L,0)e^{-i\lambda L}-e^{\lambda^2t-i\lambda l_2(t)}q(l_2(t),t)+\int_{0}^{t}e^{\lambda ^2s-i\lambda l_2(s)}q'(l_2(s),s)ds.
%\end{split}
%\label{dheat2}
%\end{eqnarray}
%\end{equation}
Using the formulae (\ref{dheat1}),
%and (\ref{dheat2}), 
equation (\ref{GRt2heat}) can be written as
\begin{eqnarray}
&e^{-i\lambda l_2(t)}q(l_2(t),t)-e^{-i\lambda l_1(t)}q(l_1(t),t)+i\lambda \int_{l_1(t)}^{l_2(t)}e^{-i\lambda x}q(x,t)dx
%\nonumber 
%\\=&
=e^{-\lambda ^2t}\Big(i\lambda \hat{q}_0(\lambda )-q(0,0)+e^{-i\lambda L}q(L,0)\Big)\nonumber 
\\&-e^{-\lambda ^2t}\int_{0}^{t}e^{-i\lambda l_1(s)+\lambda ^2s}\left[\frac{dq}{ds}(l_1(s),s)+i\lambda q_x(l_1(s),s)\right]ds
\nonumber \\&
+e^{-\lambda ^2t}\int_{0}^{t}e^{-i\lambda l_2(s)+\lambda ^2s}\left[\frac{dq}{ds}(l_2(s),s)+i\lambda q_x(l_2(s),s)\right]ds.
\label{qxrepresentationheat}
\end{eqnarray}

Using the definitions  (\ref{iv}), (\ref{bvD}), (\ref{bvN}), and the identities
\begin{subequations}
\begin{align}
&e^{-i\lambda l_2(t)}q(l_2(t),t)-e^{-i\lambda l_1(t)}q(l_1(t),t)+i\lambda
%\hat{q}(t,\lambda)
 \int_{l_1(t)}^{l_2(t)}e^{-i\lambda x}q(x,t)dx
 =\int_{l_1(t)}^{l_2(t)}e^{-i\lambda x}q_x(x,t)dx,
\\
&i\lambda \hat{q}_0(\lambda )-q(0,0)+e^{-i\lambda L}q(L,0)=\int_{0}^Le^{-i\lambda x}q'_0(x)dx,
\end{align}
\label{ids}
\end{subequations}
%and using (\ref{iv}), (\ref{bvD}), (\ref{bvN}),
equation (\ref{qxrepresentationheat}) can be written in the form
\begin{eqnarray}
\int_{l_1(t)}^{l_2(t)}e^{-i\lambda x}q_x(x,t)dx&=e^{-\lambda ^2t}\int_{0}^Le^{-i\lambda x}q'_0(x)dx-e^{-\lambda ^2t}\int_{0}^{t}e^{-i\lambda l_1(s)+\lambda ^2s}\left(f_0'(s)+i\lambda f_1(s)\right)ds
\nonumber \\&~~+e^{-\lambda ^2t}\int_{0}^{t}e^{-i\lambda l_2(s)+\lambda ^2s}\left(g_0'(s)+i\lambda g_1(s)\right)ds.
\label{qxrepresentation1heat}
\end{eqnarray}

The term on the left hand side of (\ref{qxrepresentation1heat}) is the Fourier
transform of $q_x(x, t)$ on the finite interval $l_1(t)<x<l_2(t)$.
Inverting this Fourier transform, we obtain
\begin{eqnarray}
q_x(x,t)=\frac{1}{2\pi}\int_{-\infty}^{+\infty}e^{i\lambda x-\lambda ^2t}&\Big[\int_{0}^Le^{-i\lambda \xi}q'_0(\xi)d\xi-\int_{0}^{t}e^{-i\lambda l_1(s)+\lambda ^2s}\left(f_0'(s)+i\lambda f_1(s)\right)ds
\nonumber \\&+\int_{0}^{t}e^{-i\lambda l_2(s)+\lambda ^2s}\left(g_0'(s)+i\lambda g_1(s)\right)ds\Big]d\lambda, ~~(x,t)\in \Omega(T).
\label{qxheat}
\end{eqnarray}
Assuming that $q_x(x,t) = 0$ for $x < l_1(t)$ and for $x > l_2(t)$, equation (\ref{qxheat}) is also formally valid at
$x = l_1(t)$ and at $x = l_2(t )$.
Hence, we obtain
\begin{eqnarray}
f_1(t)=\frac{1}{\pi}\int_{-\infty}^{+\infty}e^{i\lambda l_1(t)-\lambda ^2t}&\Big[\int_{0}^Le^{-i\lambda \xi}q'_0(\xi)d\xi-\int_{0}^{t}e^{-i\lambda l_1(s)+\lambda ^2s}\left(f_0'(s)+i\lambda f_1(s)\right)ds
\nonumber \\&+\int_{0}^{t}e^{-i\lambda l_2(s)+\lambda ^2s}\left(g_0'(s)+i\lambda g_1(s)\right)ds\Big]d\lambda,
\label{f1heat}
\end{eqnarray}
and
\begin{eqnarray}
g_1(t)=\frac{1}{\pi}\int_{-\infty}^{+\infty}e^{i\lambda l_2(t)-\lambda ^2t}&\Big[\int_{0}^Le^{-i\lambda \xi}q'_0(\xi)d\xi-\int_{0}^{t}e^{-i\lambda l_1(s)+\lambda ^2s}\left(f_0'(s)+i\lambda f_1(s)\right)ds
\nonumber \\&+\int_{0}^{t}e^{-i\lambda l_2(s)+\lambda ^2s}\left(g_0'(s)+i\lambda g_1(s)\right)ds\Big]d\lambda.
\label{g1heat}
\end{eqnarray}

Let
\begin{subequations}
\begin{align}
&E_j(\lambda,t,x)=e^{i\lambda (l_j(t)-x)-\lambda ^2t}, ~~j=1,2,
\label{Ejaheat}
\\
&E_{jm}(\lambda,t,s)=e^{i\lambda (l_j(t)-l_m(s))-\lambda ^2(t-s)},~~j,m=1,2.
\label{Ejbheat}
\end{align}
\label{Ej}
\end{subequations}
We rewrite equations (\ref{f1heat}) and (\ref{g1heat}) as
\begin{subequations}
\begin{align}
&\pi f_1(t)=N_1(t)-i\int_{-\infty}^{\infty}\lambda\left[\int_0^t\left(E_{11}(\lambda,t,s)f_1(s)
-E_{12}(\lambda,t,s)g_1(s)\right)ds\right]d\lambda,
\label{f1heatV}
\\
&\pi g_1(t)=N_2(t)-i\int_{-\infty}^{\infty}\lambda\left[\int_0^t\left(E_{21}(\lambda,t,s)f_1(s)
-E_{22}(\lambda,t,s)g_1(s)\right)ds\right]d\lambda,
\label{g1heatV}
\end{align}
\label{fg1heatV}
\end{subequations}
where
\begin{eqnarray}
N_j(t)=&\int_{-\infty}^{\infty}\Bigg[\int_{0}^L E_j(\lambda,t,x)q'_0(x)dx
-\int_{0}^{t}\left(E_{j1}(\lambda,t,s)f_0'(s)-E_{j2}(\lambda,t,s)g_0'(s)\right)ds\Bigg]d\lambda ,
~j=1,2.
\nonumber \\
\label{Nheat}
\end{eqnarray}

\begin{claim}
The functions $N_j(t)$, $j=1,2$, are given by (\ref{N1heat}).
\end{claim}

Interchanging the order of integration in (\ref{Nheat}), we find
\begin{eqnarray}
N_j(t)=&\int_0^L\left(\int_{-\infty}^{\infty}E_j(\lambda,t,x)d\lambda \right)q'_0(x)dx
\nonumber \\
&-\int_{0}^{t}\left[\left(\int_{-\infty}^{\infty} E_{j1}(\lambda,t,s)d\lambda\right) f_0'(s)-\left(\int_{-\infty}^{\infty} E_{j2}(\lambda,t,s)d\lambda\right) g_0'(s)\right]ds ,~~j=1,2.
\label{Nheatproof1}
\end{eqnarray}
The $\lambda$-integrals appearing in (\ref{Nheatproof1}) can be evaluated explicitly:
\begin{subequations}
\begin{eqnarray}
&&\int_{-\infty}^{\infty}E_j(\lambda,t,x)d\lambda=\frac{\sqrt{\pi}}{\sqrt{t}}e^{-\frac{(x-l_j(t))^2}{4t}}, ~~j=1,2,
\label{Ejvaluesaheat}
\\
&&\int_{-\infty}^{\infty} E_{jm}(\lambda,t,s)d\lambda=\frac{\sqrt{\pi}}{\sqrt{t-s}}e^{-\frac{(l_j(t)-l_m(s))^2}{4(t-s)}},~~j,m=1,2.
\label{Ejvaluesbheat}
\end{eqnarray}
\label{Ejvaluesheat}
\end{subequations}
Substituting the above expressions into (\ref{Nheatproof1}), we immediately obtain the formulae (\ref{N1heat}).

\begin{claim} \label{claim3.2}For a given function $h(s)\in{\bf C}[0,T]$,  the following identities hold:
\begin{eqnarray}
-i\int_{-\infty}^{\infty}\lambda \int_0^tE_{jj}(\lambda,t,s)h(s)ds
%-E_{j2}(\lambda,t,s)g_1(s)\right)ds\right]
d\lambda=\int_{0}^{t}K_{jj}(t,s)h(s)ds,
%-K_{j2}(t,s)g_1(s)\Big]ds,
\label{claim2}
\end{eqnarray}
where $K_{jj}(t,s)$, $j=1,2$, as given by (\ref{Kheat}), is a weakly singular kernel;
%\end{claim}
%
%\begin{claim} For a sufficiently smooth function $h(s)$, $s\in[0,t]$, the following identities hold:
\begin{eqnarray}
-i\int_{-\infty}^{\infty}\lambda \int_0^tE_{jm}(\lambda,t,s)h(s)ds
%-E_{j2}(\lambda,t,s)g_1(s)\right)ds\right]
d\lambda=\int_{0}^{t}K_{jm}(t,s)h(s)ds,\qquad  
%-K_{j2}(t,s)g_1(s)\Big]ds,
\label{claim3.1}
\end{eqnarray}
%where $K_{12}(t,s)$ is given by (\ref{K1heat}), and is a regular kernel.
%
%Similarly,
%\begin{eqnarray}
%\begin{split}
%-i\int_{-\infty}^{\infty}\lambda \int_0^tE_{21}(\lambda,t,s)h(s)ds
%%-E_{j2}(\lambda,t,s)g_1(s)\right)ds\right]
%d\lambda=\int_{0}^{t}K_{21}(t,s)h(s)ds,\qquad  
%%-K_{j2}(t,s)g_1(s)\Big]ds,
%\end{split}
%\label{claim3.2}
%\end{eqnarray}
where $K_{jm}(t,s)$  are the non-singular integral kernels  given by (\ref{Kheat}) when $j\neq m$.
\end{claim}

To show both claims,  we formally  interchange the order of double integration in the left hand side of (\ref{claim2}) and (\ref{claim3.1}).  We will justify the validity of this procedure by showing that all integrands are integrable.

By definition,
$$
-i\int_{-\infty}^{\infty}\lambda\int_0^tE_{jm}(\lambda,t,s)h(s)ds\,d\lambda=
\int_{-\infty}^{\infty}\lambda\int_0^t-ie^{-\lambda ^2\left(t-s\right)+i\lambda \left(l_j(t)-l_m(s)\right)}h(s)
ds\,d\lambda.
$$
We now interchange the order of integration, and use integration by parts  for the inner integral to obtain
%\begin{eqnarray}
%K_{jm}(t,s)&=&
$$
-i\int_{-\infty}^{\infty}\lambda e^{-\lambda ^2\left(t-s\right)+i\lambda \left(l_j(t)-l_m(s)\right)}
d\lambda=\frac{i}{2(t-s)}\int_{-\infty}^{\infty}  \left(\frac{\partial}{\partial\lambda}e^{-\lambda ^2(t-s)+i\lambda (l_j(t)-l_m(s))}\right)d\lambda
$$
\begin{equation}
+\frac{l_j(t)-l_m(s)}{2(t-s)}\int_{-\infty}^{\infty}  e^{-\lambda ^2(t-s)+i\lambda (l_j(t)-l_m(s))}d\lambda , \quad~~j,m=1,2.
\label{Kheatc}
\end{equation}
The first integral of the right hand side of the above equation vanishes,
whereas the second integral can be computed explicitly as (see (\ref{Ejvaluesbheat}))
$$
\frac{l_j(t)-l_m(s)}{2(t-s)}\int_{-\infty}^{\infty}  e^{-\lambda ^2(t-s)+i\lambda (l_j(t)-l_m(s))}d\lambda=
\frac{l_j(t)-l_m(s)}{2(t-s)}\frac{\sqrt{\pi}}{\sqrt{t-s}}e^{-\frac{(l_j(t)-l_m(s))^2}{4(t-s)}}, \quad ~~j,m=1,2.
$$
which is the kernel  $K_{jm}(t,s)$ as defined by equation (\ref{Kheat}).

\noindent
Note that 
\begin{itemize}
\item[]
{\bf If $j=m$:}

\noindent
In this case, the singularity at $s=t$ due to the term $\frac1{t-s}$ is removable, as 
$$
\lim_{s\to t}\frac{l_j(t)-l_j(s)}{2(t-s)}=\frac 1 2 l_j'(t).
$$
Hence the kernel $K_{jj}$ has the weak, integrable singularity $\frac1{\sqrt{t-s}}$.  
\item[]
{\bf If $j\neq m$: }

\noindent
In this case, the singularity $\frac1{(t-s)^{3/2}}$ is removable as it is cancelled by the zero of the exponential term $e^{-\frac{(l_j(t)-l_m(s))^2}{4(t-s)}}$. Therefore  the kernel $K_{jm}$, $j\neq m$, is regular at $s=t$.
\end{itemize}

\medskip
Collecting the results of these claims, equations (\ref{Nheat}) yield the system (\ref{qxrepresentationl1F}).

\begin{remark}
Our requirement in Claim \ref{claim3.2} that $h \in C[0,T]$ is sufficient to justify the interchange in the integration order, and it is all we need for our purposes, as this is the class of the unknown boundary values we seek to determine.
\end{remark}

\subsection{The linear Schr\"{o}dinger equation}

We consider the linear Schr\"{o}dinger equation (\ref{ls}). As for the case of the heat equation considered in the previous section, to 
%formulated in the time-dependent domain (\ref{MI}), with the given initial value (\ref{iv}) and Dirichlet boundary conditions (\ref{bvD}).
%The functions $f_1(t)$ and $g_1(t)$ denote the unknown Neumann boundary values at $x=l_1(t)$ and $x=l_2(t)$ respectively.
%
%We 
determine the unknown Neumann boundary values $f_1(t)$ and $g_1(t)$ given in (\ref{bvN}) in terms of the given functions $q_0(x)$, $f_0(t)$ and $g_0(t)$ we solve the global relation for the unknown boundary values $f_1(t)$ and $g_1(t)$.

In the case of linear Schr\"{o}dinger equation, we have
\begin{eqnarray}
\omega(\lambda )=\lambda ^2, ~~Q(x,t)=iq_x(x,t)-\lambda q(x,t).
\label{lsOQ}
\end{eqnarray}
Hence the global relation (\ref{GRt2}) becomes
 \begin{eqnarray}
&\hat{q}_0(\lambda )-e^{i\lambda ^2t}\hat{q}(t,\lambda )-\int_{0}^{t}e^{-i\lambda l_1(s)+i\lambda ^2s}\left[\Big(l'_1(s)-\lambda \Big)q(l_1(s),s)+iq_x(l_1(s),s)\right]ds
\nonumber \\&+\int_{0}^{t}e^{-i\lambda l_2(s)+i\lambda ^2s}\left[\Big(l'_2(s)-\lambda \Big)q(l_2(s),s)+iq_x(l_2(s),s)\right]ds=0, ~~\lambda \in \mathbb{C},~~0<t<T.
\label{GRt2ls}
\end{eqnarray}

In analogy with the case of the heat equation, in order to obtain a representation of $q_x (x, t)$, we multiply (\ref{GRt2ls})
by $i\lambda $ and then employ integration by parts for the terms involving the known functions $q(l_1(t),t)$ and $q(l_2(t),t)$.
In this way, equation (\ref{GRt2ls}) yields the equation
\begin{eqnarray}
&e^{-i\lambda l_2(t)}q(l_2(t),t)-e^{-i\lambda l_1(t)}q(l_1(t),t)+i\lambda \int_{l_1(t)}^{l_2(t)}e^{-i\lambda x}q(x,t)dx
%\nonumber \\=&
=e^{-i\lambda ^2t}\Big(i\lambda\hat{q}_0(\lambda )-q(0,0)+e^{-i\lambda L}q(L,0)\Big)
\nonumber \\&-e^{-i\lambda ^2t}\int_{0}^{t}e^{-i\lambda l_1(s)+i\lambda ^2s}\left[\frac{dq}{ds}(l_1(s),s)-\lambda q_x(l_1(s),s)\right]ds
%\nonumber \\&
+e^{-i\lambda ^2t}\int_{0}^{t}e^{-i\lambda l_2(s)+i\lambda ^2s}\left[\frac{dq}{ds}(l_2(s),s)-\lambda q_x(l_2(s),s)\right]ds.
\nonumber \\
\label{qxrepresentation}
\end{eqnarray}
Employing the identities (\ref{ids}), equation (\ref{qxrepresentation}) can be written in the form
\begin{eqnarray}
&\int_{l_1(t)}^{l_2(t)}e^{-i\lambda x}q_x(x,t)dx
%\nonumber\\=&
=e^{-i\lambda ^2t}\int_{0}^Le^{-i\lambda \xi}q'_0(\xi)d\xi-e^{-i\lambda ^2t}\int_{0}^{t}e^{-i\lambda l_1(s)+i\lambda ^2s}\left[\frac{dq}{ds}(l_1(s),s)-\lambda q_x(l_1(s),s)\right]ds
\nonumber \\&
+e^{-i\lambda ^2t}\int_{0}^{t}e^{-i\lambda l_2(s)+i\lambda ^2s}\left[\frac{dq}{ds}(l_2(s),s)-\lambda q_x(l_2(s),s)\right]ds.
\label{qxrepresentation1}
\end{eqnarray}

The term on the left hand side of (\ref{qxrepresentation1}) is the Fourier
transform of $q_x(x, t)$.
Inverting this Fourier transform of $q_x(x,t)$ and evaluating it at $x =l_j(t)$, $j=1,2$, we arrive at
\begin{eqnarray}
&q_x(l_j(t),t)
\nonumber\\=&\frac{1}{\pi}\int_{-\infty}^{\infty}e^{i\lambda l_j(t)}\Bigg[e^{-i\lambda ^2t}\int_{0}^Le^{-i\lambda x}q'_0(x)dx
-e^{-i\lambda ^2t}\int_{0}^{t}e^{-i\lambda l_1(s)+i\lambda ^2s}\Big(\frac{dq}{ds}(l_1(s),s)-\lambda q_x(l_1(s),s)\Big)ds
\nonumber\\&+e^{-i\lambda ^2t}\int_{0}^{t}e^{-i\lambda l_2(s)+i\lambda ^2s}\Big(\frac{dq}{ds}(l_2(s),s)-\lambda q_x(l_2(s),s)\Big)ds\Bigg]d\lambda, ~~j=1,2.
\label{qxrepresentationl12}
\end{eqnarray}

Setting
\begin{subequations}
\begin{align}
&\mathcal{E}_j(\lambda,t,x)=e^{i\lambda (l_j(t)-x)-i\lambda ^2t}, ~~j=1,2,
\label{Ejals}
\\
&\mathcal{E}_{jm}(\lambda,t,s)=e^{i\lambda (l_j(t)-l_m(s))-i\lambda ^2(t-s)},~~j,m=1,2,
\label{Ejbls}
\end{align}
\label{Ejls}
\end{subequations}
and using (\ref{iv}), (\ref{bvD}) and (\ref{bvN}),
we write equation (\ref{qxrepresentationl12}) in the form
\begin{subequations}
\begin{align}
&\pi f_1(t)=N_1(t)+\int_{-\infty}^{\infty}\lambda\left[\int_0^t\left(\mathcal{E}_{11}(\lambda,t,s)f_1(s)
-\mathcal{E}_{12}(\lambda,t,s)g_1(s)\right)ds\right]d\lambda,
\label{qxrepresentationlsa}
\\
&\pi g_1(t)=N_2(t)+\int_{-\infty}^{\infty}\lambda\left[\int_0^t\left(\mathcal{E}_{21}(\lambda,t,s)f_1(s)
-\mathcal{E}_{22}(\lambda,t,s)g_1(s)\right)ds\right]d\lambda,
\label{qxrepresentationlsb}
\end{align}
\label{qxrepresentationls}
\end{subequations}
where, for $j=1,2$,
\begin{eqnarray}
N_j(t)=\int_{-\infty}^{\infty}\Bigg[\int_{0}^L\mathcal{E}_j(\lambda,t,x)q'_0(x)dx
-\int_{0}^{t}\left(\mathcal{E}_{j1}(\lambda,t,x)f_0'(s)-\mathcal{E}_{j2}(\lambda,t,x)g_0'(s)\right)ds\Bigg]d\lambda.
\label{Nls}
\end{eqnarray}
We want to rewrite (\ref{qxrepresentationls}) in the form  (\ref{qxrepresentationls1F}).
\begin{claim}
The functions $N_j(t)$, $j=1,2$, are given by  (\ref{N1ls}).
%\begin{eqnarray}
%N_j(t)=\frac{(1-i)\sqrt{2\pi}}{2}\Bigg\{\frac{1}{\sqrt{t}}\int_{0}^{L}e^{\frac{i(x-l_j(t))^2}{4t}}q'_0(x)dx
%-\int_{0}^{t}\Big[\frac{e^{\frac{i(l_j(t)-l_1(s))^2}{4(t-s)}}}{\sqrt{t-s}}f_0'(s)
%-\frac{e^{\frac{i(l_j(t)-l_2(s))^2}{4(t-s)}}}{\sqrt{t-s}}g_0'(s)\Big]ds\Bigg\}.
%\label{N1ls}
%\end{eqnarray}
\end{claim}

Interchanging the order of integration in (\ref{Nls}), we find
\begin{eqnarray}
N_j(t)&=&\int_0^L\left(\int_{-\infty}^{\infty}\mathcal{E}_j(\lambda,t,x)d\lambda \right)q'_0(x)dx
\nonumber\\
&-&\int_{0}^{t}\left[\left(\int_{-\infty}^{\infty} \mathcal{E}_{j1}(\lambda,t,s)d\lambda\right) f_0'(s)-\left(\int_{-\infty}^{\infty} \mathcal{E}_{j2}(\lambda,t,s)d\lambda\right) g_0'(s)\right]ds ,~~j=1,2.
\label{Nlsproof1}
\end{eqnarray}
The $\lambda$-integrals appearing in (\ref{Nlsproof1}) can be evaluated explicitly:
\begin{subequations}
\begin{eqnarray}
&&\int_{-\infty}^{\infty}\mathcal{E}_j(\lambda,t,x)d\lambda=\frac{(1-i)\sqrt{2\pi}}{2\sqrt{t}}e^{\frac{i(x-l_j(t))^2}{4t}}, ~~j=1,2,
\label{Ejlsvaluesa}
\\
&&\int_{-\infty}^{\infty} \mathcal{E}_{jm}(\lambda,t,s)d\lambda=\frac{(1-i)\sqrt{2\pi}}{2\sqrt{t-s}}e^{\frac{i(l_j(t)-l_m(s))^2}{4(t-s)}},~~j,m=1,2.
\label{Ejlsvaluesb}
\end{eqnarray}
\label{Ejlsvalues}
\end{subequations}
Substituting (\ref{Ejlsvalues}) into (\ref{Nlsproof1}) we immediately obtain the formulae (\ref{N1ls}).

\begin{claim}  For a given function $h(s)\in{\bf C}[0,T]$, the following identities hold:
\begin{eqnarray}
\int_{-\infty}^{\infty}\lambda\int_0^t\mathcal{E}_{jj}(\lambda,t,s)h(s)
%-\mathcal{E}_{j2}(\lambda,t,s)g_1(s)\right)ds\right]
dsd\lambda=\int_{0}^{t}K_{jj}(t,s)h(s)
%-K_{j2}(t,s)g_1(s)\Big]
ds,
\label{claim4}
\end{eqnarray}
\begin{eqnarray}
\int_{-\infty}^{\infty}\lambda\int_0^t\mathcal{E}_{jm}(\lambda,t,s)h(s)
%-\mathcal{E}_{j2}(\lambda,t,s)g_1(s)\right)ds\right]
dsd\lambda=\lim_{\varepsilon\to 0}\int_{0}^{t}K_{jm}(t,s,\varepsilon)h(s)
%-K_{j2}(t,s)g_1(s)\Big]
ds,
\label{claim5}
\end{eqnarray}
where $K_{jj}(t,s)$, $K_{jm}(s,t,\varepsilon)$ $j,m=1,2$, are the integral kernels given by (\ref{K1ls}) and (\ref{K2ls}) respectively.
%\begin{eqnarray}
%\begin{split}
%K_{jj}(t,s)=\frac{(1-i)\sqrt{2\pi}}{4}\frac{l_j(t)-l_j(s)}{t-s}
%\frac{e^{\frac{i(l_j(t)-l_j(s))^2}{4(t-s)}}}{\sqrt{t-s}},
%~~\varepsilon>0,~~0<s<t<T, \end{split}
%\label{K1ls}
%\end{eqnarray}
%and $K_{jm}(t,s,\varepsilon)$, $j,m=1,2$, $j\neq m$ is given by
%\begin{eqnarray}
%\begin{split}
%K_{jm}(t,s,\varepsilon)=\frac{(1-i)\sqrt{2\pi}}{4}\frac{l_j(t)-l_m(s)}{t-s-i\varepsilon}
%\frac{e^{\frac{i(l_j(t)-l_m(s))^2}{4(t-s-i\varepsilon)}}}{\sqrt{t-s-i\varepsilon}},
%~~\varepsilon>0,~~0<s<t<T, \end{split}
%\label{K2ls}
%\end{eqnarray}
\end{claim}

If we interchange the order of double integration in the left hand side of (\ref{claim4}),  we obtain an integrable functions of $s$. However, for the case $j\neq m$ of (\ref{claim5}),  we obtain a function that is not integrable with respect to $s$. Thus, before interchanging the order of the integration we must first regularise the relevant $\lambda$-integral. Therefore, we write
\begin{eqnarray}
&\int_{-\infty}^{\infty}\int_0^t\lambda\mathcal{E}_{jm}(\lambda,t,s)h(s)d s
%-\mathcal{E}_{j2}(\lambda,t,s)g_1(s)\right)ds\right]
d\lambda
=\int_{-\infty}^{\infty}\lambda \int_0^t\lim_{\varepsilon\rightarrow 0^{+}}e^{-\varepsilon\lambda ^2-i\lambda ^2\left(t-s\right)+i\lambda \left(l_j(t)-l_1(s)\right)}h(s)ds
%-e^{-\varepsilon\lambda ^2-i\lambda ^2\left(t-s\right)+i\lambda \left(l_j(t)-l_2(s)\right)}g_1(s)\right)ds\right]
d\lambda
\nonumber\\
=&\lim_{\varepsilon\rightarrow 0^{+}}\int_{-\infty}^{\infty}\lambda \int_0^t e^{-\varepsilon\lambda ^2-i\lambda ^2\left(t-s\right)+i\lambda \left(l_j(t)-l_1(s)\right)}h(s)ds\,d\lambda,
\label{claim4F}
\label{qxrepresentationls1Fagain}
\end{eqnarray}
where the last identity follows from the dominated convergence theorem, thanks to the exponential decay of the term $e^{-\varepsilon\lambda ^2}$. Now we can interchange the order of integration, hence the expression in  (\ref{claim4F}) is equal to 
%$$
\begin{eqnarray}
\lim_{\varepsilon\rightarrow 0^{+}}
\int_{0}^{t}K_{jm}(t,s,\varepsilon)h(s)ds,\quad
%$$
%where
%\begin{eqnarray}
K_{jm}(t,s,\varepsilon)=\int_{-\infty}^{\infty}\lambda  e^{-i\lambda ^2\left(t-s-i\varepsilon\right)+i\lambda \left(l_j(t)-l_m(s)\right)}d\lambda , 
\nonumber\\
~~j,m=1,2,\,j\neq m.
\label{K}
\end{eqnarray}
The $\lambda$-integral (\ref{K}) can be evaluated explicitly.
Indeed, the $\lambda$-derivative of the exponent of the exponential appearing in (\ref{K}) is given by the expression:
$$\left[-2i\lambda(t-s-i\varepsilon)+i\left(l_j(t)-l_m(s)\right)\right]e^{-i\lambda ^2\left(t-s-i\varepsilon\right)+i\lambda \left(l_j(t)-l_m(s)\right)}.$$
Hence, we can rewrite $K_{jm}(t,s,\varepsilon)$ in the form
\begin{eqnarray}
K_{jm}(t,s,\varepsilon)=&\frac{-1}{2i(t-s-i\varepsilon)}\int_{-\infty}^{\infty}  \left(\frac{\partial}{\partial\lambda}e^{-i\lambda ^2\left(t-s-i\varepsilon\right)+i\lambda \left(l_j(t)-l_m(s)\right)}\right)d\lambda
\nonumber\\
&+\frac{l_j(t)-l_m(s)}{2(t-s-i\varepsilon)}\int_{-\infty}^{\infty}  e^{-i\lambda ^2\left(t-s-i\varepsilon\right)+i\lambda \left(l_j(t)-l_m(s)\right)}d\lambda , ~~j,m=1,2.
\label{Kc}
\end{eqnarray}
The first integral of the right hand side of the above equation vanishes, because of the large $\lambda$ decay of the term $e^{-\lambda^2\varepsilon}$, 
whereas the second integral can be computed explicitly (see (\ref{Ejlsvaluesb}).
This yields the expression
\begin{eqnarray}
K_{jm}(t,s,\varepsilon)=\frac{l_j(t)-l_m(s)}{2(t-s-i\varepsilon)}\frac{(1-i)\sqrt{2\pi}}{2\sqrt{t-s-i\varepsilon}}e^{\frac{i(l_j(t)-l_m(s))^2}{4(t-s-i\varepsilon)}}, ~~j,m=1,2, \; j\neq m.
\label{Kc1}
\end{eqnarray}

In summary, using (\ref{qxrepresentationls1Fagain})-(\ref{Kc}), we find
\begin{equation}
\int_{-\infty}^{\infty}\lambda\mathcal{E}_{jm}(\lambda,t,s)h(s)d s=\frac{(1-i)\sqrt{2\pi}}{4}\lim_{\varepsilon\rightarrow 0^{+}}
\int_{0}^{t}\frac{l_j(t)-l_m(s)}{(t-s-i\varepsilon)^{3/2}}e^{\frac{i(l_j(t)-l_m(s))^2}{4(t-s-i\varepsilon)}}h(s)ds,
\label{final}
\end{equation}
yielding the conclusion of Proposition \ref{theoremLS}.

\section{Proof of the main theorems}
We now sketch the step needed to prove that the representation derived in the previous section yields Volterra integral equations that admit a unique solution.  The only new ingredient in this section is the analysis of the $\varepsilon \to 0$ limit of the integrals appearing in the representation given by Proposition \ref{theoremLS}.

After proving that these limits yield a well defined system of Volterra integral equations, possibly weakly singular, the proof is analogous to the proof given in  \cite{AP3} for the problem formulated on $l(t)<x<\infty$. We refer to these papers for details, and concentrate on showing that the integral equations derived in the previous section  are of a type that can be treated using classical results. 

We note that for both the case of the heat equation and the linear Schr\"odinger equation, the case  that 
 $j=m$ yields a weakly singular kernel. 

Indeed, in this case, the singularity at $s=t$ due to the term $\frac1{t-s}$ is removable, as 
$$
\lim_{s\to t}\frac{l_j(t)-l_j(s)}{2(t-s)}=\frac 1 2 l_j'(t).
$$
Hence the kernel $K_{jj}$ has the weak, integrable singularity $\frac1{\sqrt{t-s}}$.  

\smallskip
However, if $j\neq m$, this is not the case. We consider the two theorems separately.

\subsection{Theorem \ref{main}}

For the kernel given by expression (\ref{Kheat}), the singularity $\frac1{(t-s)^{3/2}}$ is removable as it is cancelled by the zero of the exponential term $e^{-\frac{(l_j(t)-l_m(s))^2}{4(t-s)}}$. Therefore  the kernel $K_{jm}$, $j\neq m$, is regular at $s=t$.

Under our regularity assumptions on the known data, it follows that the system of Volterra integral equations (\ref{qxrepresentationl1F}) admits a unique solution. The proof is identical to the proof given in \cite{AP3} and relies on general results for Volterra integral equations with weakly singular kernels, given e.g. in \cite{MF}. 

\subsection{Theorem \ref{main2}}

Next we consider the kernel given by expression (\ref{Kc1}).
For $\varepsilon>0$, this kernel has no singularity, hence invoking again the general results of \cite{AP3, MF} we can deduce that the vector Volterra integral equation 
%\begin{subequations}
%\begin{align}
%&\pi f_1(t)=N_1(t)+\int_{0}^{t}K_{11}(t,s)f_1(s)ds-\int_0^tK_{12}(t,s,\varepsilon)g_1(s)ds,~~0<t<T,
%\label{1noeps}
%\\
%&\pi g_1(t)=N_2(t)+\int_{0}^{t}K_{21}(t,s,\varepsilon)f_1(s)ds-\int_0^tK_{22}(t,s)g_1(s)ds,~~0<t<T,
%\label{2noeps}
%\end{align}
%\label{noeps}
%\end{subequations}
\begin{equation}
%{\mathbf h}^\varepsilon(t)
\left(\begin{array}{c}f_1^\varepsilon(t) \\g_1^\varepsilon(t)\end{array}\right)=\left(\begin{array}{c}N_1(t)\\N_2(t)\end{array}\right)+\int_{0}^{t}
\left(\begin{array}{cc}
K_{11}(t,s)&-K_{12}(t,s,\varepsilon)
\\
K_{21}(t,s,\varepsilon)&-K_{22}(t,s)
\end{array}\right)\left(\begin{array}{c}f_1^\varepsilon(s) \\g_1^\varepsilon(s)\end{array}\right)ds
%{\mathbf h}^\varepsilon(s)ds
\label{svie}
\end{equation}
admits a unique solution $(f_1^\varepsilon, g_1^\varepsilon)\in {\bf C}[0,T)\times {\bf C}[0,T)$ for every $\varepsilon>0$.

The last step in the proof is the consideration of  the $\varepsilon \to 0^+$ limit in the expression above. It must be shown that this limit exists. To avoid technicalities and focus on the essential issue of the $\varepsilon \to 0$ limit, we show this for the analogous scalar case - the extension to the case of the vector integral equation (\ref{svie}) is immediate.

We first consider the $\varepsilon$-dependent kernel $K_{12}$. 
Recall that, for a function $h(s):[0,t]\to {\mathbb R}$ bounded and sufficiently regular, 
\begin{equation}
 \int_0^t K_{12}(t,s,\varepsilon)h(s)ds=\frac{(1-i)\sqrt{2\pi}}{4} \int_{0}^{t}\frac{l_1(t)-l_2(s)}{(t-s-i\varepsilon)^{3/2}}e^{\frac{i(l_1(t)-l_2(s))^2}{4(t-s-i\varepsilon)}}h(s)ds.
 \label{partway}\end{equation}
 We now consider the exponential appearing in the integrand
 $$
{\cal E}_{12}(t,s,\varepsilon)= e^{\frac{i(l_1(t)-l_2(s))^2}{4(t-s-i\varepsilon)}}
 $$
 Differentiating ${\cal E}_{12}$ with respect to $s$, and rearranging, we can write the integrand in (\ref{partway}) as
 $$
 \frac{l_1(t)-l_2(s)}{(t-s-i\varepsilon)^{3/2}}{\cal E}_{12}(t,s,\varepsilon)=-4i\frac{(t-s-i\varepsilon)^{1/2}}{l_1(t)-l_2(s)-2l_2'(s)(t-s-i\varepsilon)}\frac{\partial {\cal E}_{12}}{\partial s}.
 $$
 Hence
 $$
   \int_0^t K_{12}(t,s,\varepsilon)h(s)ds=-(i+1)\sqrt{2\pi}\int_0^t\frac{(t-s-i\varepsilon)^{1/2}}{l_1(t)-l_2(s)-2l_2'(s)(t-s-i\varepsilon)}\frac{\partial {\cal E}_{12}}{\partial s}h(s)ds.
  $$
 Integration by parts yields
% \begin{eqnarray}
$$
-\frac 1 {(i+1)\sqrt{2\pi}}  \int_0^t K_{12}(t,s,\varepsilon)h(s)ds=
$$
%%%%%%%%%%%%%%%%%
%%%%%%%%HERE
%%%%%%%%%%%%%%%%%
$$
\frac{(-i\varepsilon)^{1/2}}{l_1(t)-l_2(t)+2l_2'(t)i\varepsilon}h(t){\cal E}_{12}(t,t,\varepsilon)-\frac{(t-i\varepsilon)^{1/2}}{l_1(t)-l_2(0)-2l_2'(0)(t-i\varepsilon)}h(0){\cal E}_{12}(t,0,\varepsilon)
$$
$$
-\int_0^t{\cal E}_{12}(t,s,\varepsilon)\left\{\frac{-\frac{h(s)}{2(t-s-i\varepsilon)^{1/2}}+(t-s-i\varepsilon)^{1/2}h'(s)}{l_1(t)-l_2(s)-2l_2'(s)(t-s-i\varepsilon)}-\frac{(t-s-i\varepsilon)^{1/2}\left[l_2'(s)-2l_2''(s)(t-s-i\varepsilon)\right]h(s)}{\left(l_1(t)-l_2(s)-2l_2'(s)(t-s-i\varepsilon)\right)^2}
\right\}ds.
 $$
 We now need to take the limit as $\varepsilon\to 0$; in order to pass to the limit inside the integral on the right hand side using the dominated convergence theorem, we must show that the integrand is dominated by an integrable function. 
 
 Let
 \begin{equation}
H_1(t,s)=l_1(t)-l_2(s)-2l_2'(s)(t-s).
\label{hts}
\end{equation}
If $H_1(t,s)\neq0$ for all $s\in[0,t]$, then the integrand can be dominated by $g(s)=\frac {c}{\sqrt{t-s}}$,  and this function is integrable in $[0,t]$.  Therefore, under this assumption, that we will return to below, we can pass to limit under the integral  and  we find
$$\lim_{\varepsilon\to 0}\frac 1 {(i+1)\sqrt{2\pi}}  \int_0^t K_{12}(t,s,\varepsilon)h(s)ds=
\frac{t^{1/2}}{H_1(t,0)}h(0){\cal E}_{12}(t,0,0)+\int_0^t{\cal E}_{12}(t,s,0)\frac{(t-s)^{1/2}}{H_1(t,s)}h'(s)ds
$$
$$-\int_0^t{\cal E}_{12}(t,s,0)\left\{\frac{1}{2(t-s)^{1/2}H_1(t,s)}+\frac{\left[l_2'(s)-2l_2''(s)(t-s)\right](t-s)^{1/2}}{H_1(t,s)^2}\right\}
h(s)ds.
$$
Thus this limit has the form
$$
\lim_{\varepsilon\to 0}\frac 1 {(i+1)\sqrt{2\pi}}  \int_0^t K_{12}(t,s,\varepsilon)h(s)ds= M^h_{12}(t)-\int_0^t
{\cal E}_{12}(t,s)\frac 1{H_1(t,s)}\left[\frac 1 {2(t-s)^{1/2}}+\frac {\left[l_2'(s)-2l''_2(s)(t-s)\right](t-s)^{1/2}}{H_1(t,s)}\right]h(s)ds
$$
$$
+\int_0^t{\cal E}_{12}(t,s)(t-s)^{1/2}\frac 1{H_1(t,s)}h'(s)ds
$$
with 
$$M^h_{12}(t)={\cal E}_{12}(t,0,0)t^{1/2}\frac 1{H_1(t,0)}h(0)
$$
and 
$
H_1(t,s)$ given by (\ref{hts}).
%
%$$\lim_{\varepsilon\to 0}\frac 1 {(i+1)\sqrt{2\pi}}  \int_0^t K_{12}(t,s,\varepsilon)h(s)ds= M_{12}(t)+\int_0^t
%H_1(t,s)h(s)ds+\int_0^t
%H_2(t,s)h'(s)ds.
%$$
Hence for $f_1(t)=\lim_{\varepsilon\to 0} f_1^{\varepsilon}(t)$,  $g_1(t)=\lim_{\varepsilon\to 0} g_1^{\varepsilon}(t)$, using  equation (\ref{svie}) we find the integral equation
\begin{eqnarray}
f_1(t)&=&N_1(t)-
(i+1)\sqrt{2\pi}M^{g_1}_{12}(t)+\int_0^tK_{11}(t,s)f_1(s)ds
\nonumber \\
&+&(i+1)\sqrt{2\pi}\left\{\int_0^t
{\cal E}_{12}(t,s)\frac 1{H_1(t,s)}\left[\frac 1 {2(t-s)^{1/2}}+(t-s)^{1/2}\frac {\left[l_2'(s)-2l''_2(s)(t-s)\right](t-s)^{1/2}}{H_1(t,s)}\right]g_1(s)ds\right.
\nonumber \\
&-&\left.\int_0^t{\cal E}_{12}(t,s)(t-s)^{1/2}\frac 1{H_1(t,s)}g_1'(s)ds\right\},
\label{vls12}
\end{eqnarray}
where $N_1(t)$ is given by (\ref{N1ls}).

An analogous computation for the kernel $K_{21}$ yields
\begin{eqnarray}
g_1(t)&=&N_2(t)-
(i+1)\sqrt{2\pi}M_{21}^{f_1}(t)-\int_0^tK_{22}(t,s)g_1(s)ds
\nonumber \\
&-&(i+1)\sqrt{2\pi}\left\{\int_0^t
{\cal E}_{21}(t,s)\frac 1{H_2(t,s)}\left[\frac 1 {2(t-s)^{1/2}}+\frac {\left[l_1'(s)-2l''_1(s)(t-s)\right](t-s)^{1/2}}{H_2(t,s)}\right]f_1(s)ds\right.
\nonumber \\
&-&\left.\int_0^t{\cal E}_{21}(t,s)(t-s)^{1/2}\frac 1{H_2(t,s)}f'_1(s)ds\right\},
\label{vls21}
\end{eqnarray}
with $N_2(t)$ given by (\ref{N1ls}),
$$
{\cal E}_{21}(t,s,\varepsilon)= e^{\frac{i(l_2(t)-l_1(s))^2}{4(t-s-i\varepsilon)}},
 \qquad M^h_{21}(t)={\cal E}_{21}(t,0,0)t^{1/2}\frac 1{H_2(t,0)}h(0)
$$
 and
 $$
 H_2(t,s)=l_2(t)-l_1(s)-2l_1'(s)(t-s).$$
 
We claim that the two equations (\ref{vls12})-(\ref{vls21}) above are a system of a generalised Volterra integral equation of the second kind  with a weakly integral kernel.  

We first note that these equations are not in the usual form of a Volterra integral equation for the functions $f(t)$, $g(t)$, since the right hand side contains not only the functions but also their first derivative.  A modification of the iterative proof of existence of a solution for the usual Volterra case also works in this generalised case, see \cite{brunner}.

It remains to prove that the kernels appearing in the two integral on the right hand side of (\ref{vls12})-(\ref{vls21})  are weakly singular. This is clearly the case provided $H_1(t,s)$ and $H_2(t,s)$  do not vanish for any $s\in[0,t]$.

For $H_1(t,s)$ the condition is that
$$
l_1(t)-l_2(s)- 2l_2'(s)(t-s)\neq 0,\quad  \forall \,0\leq s\leq t.
$$
Since $l_1(t)<l_2(t)$, if we assume that $l_2(t)$ satisfies the condition (\ref{bcurves}) we have
$$
l_1(t)-l_2(s)- 2l_2'(s)(t-s)<l_2(t)-l_2(s)-2l'(s)(t-s)=(t-s)[l_2'(\sigma)-2l'(s)]<0,\quad s\leq \sigma\leq t.
$$
Similarly, for $H_2(t,s)$ the condition is
$$
l_2(t)-l_1(s)- 2l_1'(s)(t-s)\neq 0,\quad  \forall \,0\leq s\leq t.
$$
and if $l_1(t)$ satisfies the condition (\ref{bcurves}) we have
$$
l_2(t)-l_1(s)- 2l_1'(s)(t-s)>l_1(t)-l_1(s)-2l_1'(s)(t-s)>0.
$$
Hence under the assumption (\ref{bcurves}),  the regularity  condition is satisfied for both $K_{12}$ and $K_{21}$. Hence the assumption  (\ref{bcurves}) is sufficient to ensure the sought regularity. This  completes the proof of the theorem.

\begin{remark}[Linear boundaries]

The case when boundaries are linear is often of interest in applications. In this case, we can prove the main result {\em without} the need to assume equations (\ref{bcurves}).This shows that (\ref{bcurves}) is a sufficient but not necessary condition. 

Indeed, assume that the  boundary curves are of the form
\begin{equation}
l_1(s)=\alpha s, \qquad l_2(s)=\beta s+L, \qquad 0<s<t,\;\; 0<2\alpha<\beta,\;\; L\geq 0.
\label{lines}
\end{equation}
For $H_1(t,s)$ the non-zero condition becomes 
$$
(\alpha-2\beta)t+\beta s-L\neq 0 \Longleftrightarrow s\neq \frac {2\beta-\alpha}{\beta}t+\frac L\beta.
$$
Since $s\leq t$, this  always holds if $\beta>\alpha>0$. 
For $H_2(t,s)$, under the assumption (\ref{lines}), 
the non-zero condition is  
$$
(\beta-2\alpha)t+\alpha s+L\neq 0
$$
which is always true if all terms are positive, i.e. if $\beta>2\alpha>0$.

Therefore, Theorem (\ref{main2}) is valid for linear boundaries of the form (\ref{lines}).
\end{remark}

\section{Conclusions}
We have shown how to give a solution representation for boundary value problems for linear evolution equation in one spatial variable, posed between two time-dependent boundaries, issuing from the common point set at the origin of the $(x,t)$ plane. 

For the specific example of the heat equation, the solution is obtained as the unique solution of a system of Volterra integral equations (\ref{qxrepresentationl1F}), valid for any choice of differentiable  boundary curves not intersecting for positive times. 

For the case of the linear Schr\"odinger equation, the solution is again given as the unique solution of the system of generalised Volterra integral equations (\ref{vls12})-(\ref{vls21}), but only for a more restricted class of boundaries, satisfying condition (\ref{bcurves}) or (\ref{lines}).

\section*{ACKNOWLEDGMENTS}

A.S. Fokas was supported by EPSRC in the form of a senior fellowship.
B. Pelloni was supported by EPSRC grant F15R10379.
B. Xia was supported by the National Natural Science Foundation of China, Grant No. 11771186.

\small{

}
\end{document}